\documentstyle [12pt]{article}
\newenvironment{mproof}{\begin{trivlist}\item[]{\em
Proof: }}{\hfill$\Box$\end{trivlist}}
\def \IR{\hbox{{\rm I}\kern-.2em\hbox{{\rm R}}}}
\def \iR{\hbox{{\sevenrm I\kern-.2em\hbox{\sevenrm R}}}}
\def \IN{\hbox{{\rm I}\kern-.2em\hbox{\rm N}}}
\def \IC{\hbox{{\rm I}\kern-.6em\hbox{\bf C}}}
\def \IQ{\hbox{{\rm I}\kern-.6em\hbox{\bf Q}}}
\def \ZZ{\hbox{{\rm Z}\kern-.4em\hbox{\rm Z}}}

\begin{document}
\newcommand{\STAR}{*}
\newcommand{\dbq}{Q}
\newcommand{\Z }{Z}
\newcommand{\ca}{\bigcap}
\newcommand{\rarr}{\longrightarrow}
\newcommand{\ib}{\subseteq}
\newcommand{\ap}{\rightarrow}
\newcommand{\cx}{\bigotimes}
\newcommand{\memof}{\in}
\newcommand{\cir}{\circ}

\begin{center}
\bf{Algebraic orders on $K_0 $ and approximately finite \\
operator algebras \footnote{REVISED DECEMBER 1992}}\\
\end{center}
\begin{center}
Stephen C. Power\\
Department of Mathematics \\
University of Lancaster \\
England LA1 4YF \\

\end{center}
\par\vspace{2.0\baselineskip}
\par
Approximately finite (AF) $C^{\STAR} $-algebras are
classified by approximately finite ($r$-discrete principal) groupoids.
Certain natural triangular subalgebras of AF $C^{\STAR} $-algebras are
similarly classified by triangular subsemigroupoids of AF groupoids [10].
Putting this in a more intuitive way, such subalgebras $A$ are
classified by the topologised fundamental binary relation
$R(A)$ induced on the Gelfand space of the
masa $A \ \ca\	A^{\STAR} $ by the normaliser of
$A \ \ca\  A^{\STAR} $ in $A$.	(This relation $R(A)$ is also determined by
any matrix unit system for $A$ affiliated with $A \ \ca\  A^{\STAR} $.)
The fundamental relation $R(A)$ has been useful both in understanding the
isomorphism classes of specific algebras and in the general
structure theory of triangular and
chordal subalgebras of AF $C^{\STAR} $-algebras
(\cite{h+p}, \cite{p-triangular} ,\cite{p-tensor},\cite{th-triangular},
\cite{th-dilation}).	Nevertheless it is desirable to have more
convenient and computable invariants associated with the
$K_0 $ group, and we begin such an inquiry in this paper.
\par
We introduce the algebraic order and the strong algebraic
order on the scale of the $K_0 $ group
of a non-self-adjoint subalgebra of a $C^{\STAR} $-algebra.  Analogues
(and generalisations) of Elliott's classification of AF $C^{\STAR} $-algebras
are obtained for limit algebras of direct systems
$$
A_1 \rarr A_2 \rarr \ \it\hbox{...}
$$
of finite-dimensional CSL algebras (poset algebras) with
respect to certain embeddings with $C^{\STAR} $-extensions which, in
a certain sense, preserve the algebraic order.	We also
require that the systems have a certain conjugacy property.  Despite
the restrictions there are many interesting applications.  For
example conjugacy properties prevail for certain embeddings
of finite-dimensional nest algebras (block upper triangular matrix
algebras) and for systems associated with
ordered Bratteli diagrams.
\par
Hitherto the study of non-self-adjoint subalgebras of AF
$C^{\STAR} $-algebras has focused on triangular subalgebras, where
$A \ \ca\  A^{\STAR} $ is a certain
approximately finite regular maximal abelian self-adjoint algebra
(\cite{bak}, \cite{h+p}, \cite{ppw}, \cite{p+w1}, \cite{p-tensor},
\cite{p-triangular}).  See also \cite{m+s} .  However from
the point of view of identifying the algebraically
ordered scaled ordered dimension group, the viewpoint of
this paper, it is the
\it nontriangular
\rm subalgebras
which are particularly interesting since in this case
$K_0 (A)$ (which agrees with $K_0 (A \ \ca\  A^{\STAR} )$)
can be a `small group', such as $ \ZZ ^5 $ or
$ \IQ^2 $.  In such settings, the algebraic orders can be
revealed more explicitly.  For example, in Example 4.5
we have the situation in which $A \ \ca\  A^{\STAR} $ is a
simple $C^{\STAR} $-algebra in the simple $C^{\STAR} $-algebra
$C^{\STAR} (A)$, and $A$ is one of only finitely many algebras
between $A \ \ca\  A^{\STAR} $ and $C^{\STAR} (A)$.
The algebraic order of such an algebra  corresponds to partial orders
on the fibres of the surjection
$i_{\STAR} \,:\, K_0 (A) \rarr K_0 (C^{\STAR} (A))$.
\par
In section 1 we define the reflexive transitive antisymmetric
order $S(A)$ on the scale of the $K_0 $ group of a subalgebra
of a $C^{\STAR} $-algebra, and we recall some basic facts
concerning (regular) canonical subalgebras of AF $C^{\STAR} $-algebras.
In section 2 we discuss various kinds of embeddings of
finite-dimensional algebras of matrices, and we define the
strong algebraic order $S_1 (A)$ associated with a
canonical masa.  In section 3 we obtain the main results,
Theorems 3.1 and 3.2, together with various
examples and associated remarks.  In particular we consider a class of
triangular algebras associated with
\it ordered
\rm Bratteli diagrams.	In
section 4 we consider examples of non-self-adjoint subalgebras of
AF $C^{\STAR} $-algebras with small
$K_0 $ group.  In this connection we look
at
\it stationary pairs
\rm of AF $C^{\STAR} $-algebras
$D \subseteq B$. This situation includes
the context of the example mentioned above.
\par
Most of this paper was completed during a visit to the University of
Waterloo in 1989, and the author would like to thank Ken Davidson for
some useful discussions. The present preprint replaces an earlier one
which contained a mistaken (nontriangular) version of Lemma 2.5.


\noindent \bf 1. Algebra orders on $K_0$. \rm

\par\vspace{1.0\baselineskip}
\par
We start by recalling some terminology and properties of subalgebras of AF
$C^{\STAR} $-algebras.
\par
A finite-dimensional commutative subspace lattice algebra $A$,
or FDCSL algebra, is an operator algebra on a finite-dimensional
Hilbert space which contains a maximal abelian self-adjoint algebra (masa).
We say that a masa $C$ in an AF $C^{\STAR} $-algebra $B$ is
a \it canonical
\rm masa (or, more precisely, a regular canonical masa) if there is a
chain of finite-dimensional $C^{\STAR} $-subalgebras
$B_1 \ \ib\  B_2 \ \ib \ \it\hbox{...} $, with dense union, such that
the algebras
$C_k = B_k \cap C$ are masas in the algebras
$B_k $,  with dense union in $C$,
and such that, for each  $k$, the normaliser of $C_k$ in $B_k$ is contained
in the normaliser of $C_{k+1} $ in $B_{k+1} $.
A closed subalgebra $A$ of $B$ is  said to be a
\it (regular) canonical
\rm subalgebra if
$C \ \ib\  A \ \ib\  B$ for some canonical masa $C$.
In this case $A$ is necessarily the
closed union of the FDCSL algebras
$A_n \ =\  B_n \ \ca\  A$. (See \cite{p-book}.)	In particular the algebra $A$
is an approximately finite operator algebra and
is identifiable with the direct limit Banach algebra
$ {\displaystyle\lim_\to} \, A_n $ where the embeddings
possess star extensions.  Of course the converse
is true ; if $A_1 \rarr A_2 \rarr \ \it\hbox{...}$ is a direct
system of FDCSL algebras with respect to embeddings, not necessarily
unital, which have star extensions
$C^{\STAR} (A_k )\ \ \to
\ \ C^{\STAR} (A_{k+1} )$, then the Banach algebra
$A \ =\  {\displaystyle \lim_\to} \, A_k $ is completely isometrically
isomorphic to a subalgebra of the AF $C^{\STAR} $-algebra
$B \ =\  {\displaystyle \lim_\to} \, C^{\STAR} (A_k )$.
(However such a subalgebra need not be a regular canonical subalgebra
in the sense above.)
\par
The masas above coincide with approximately finite Cartan
subalgebras of AF $C^{\STAR} $-algebras \cite{ren}.
A useful discussion of them is given in the notes of
Stratila and Voiculescu \cite{s+v}.
\par
We now give a definition of $K_0 (A)$ for a not
necessarily self-adjoint subalgebra $A$ of a $C^{\STAR} $-algebra.
In the unital, or stably unital case, $K_0 (A)$ coincides with
the usual definition in terms of the stable algebraic equivalence of
idempotents.  (See Proposition 5.5.5 of [2].)  Write
$p \ap q$, or $p \ap_A q$, or $p \ap_v q$,
for $p,q$ in $ {\rm Proj} \, (A)$, the set of self-adjoint
projections of $A$, if there exists a
partial isometry $v$ in $A$ with
$v^{\STAR} v \ =\  p$, $vv^{\STAR} \ =\  q$.
Write $p \sim q$ if $v$ can be chosen in $A \ \ca\  A^{\STAR} $.
Define $K_0^+ (A)$ as the set of (Murray von Neumann)
equivalence classes $[p]$ of projections
in $ {\rm Proj}\, (A \ \cx\  M_n )$,
$n=1,2,... .$, with the usual identifications
 $A\  \cx\ M_n \subseteq A \cx M_{n+1} ,\ n =\ 1,\ 2,...\ $
A semigroup operation is given by
$[p] + [q] \ =\  [p+q]$ ,where $p$ and $q$
are representatives with $pq \ =\  0$, and
$K_0 (A)$ is, by definition, the
Grothendieck group of $K_0^+ (A)$.  For a
canonical AF subalgebra $K_0^+ (A)$ has cancellation and
embeds injectively in $K_0 (A)$.
The scale of $A$ in $K_0 (A)$ is the partially ordered
set
$ \Sigma (A) \ =\  \{[p] : p \memof {\rm Proj} (A)\} $.  A
celebrated theorem
of G. Elliott \cite{ell}  asserts that AF $C^{\STAR} $-algebras
$B_1 $ and $B_2 $ are isomorphic if there is a group
isomorphism $ \theta \,:\, K_0 (B_1 ) \rarr K_0 (B_2 )$
with $ \theta \, ( \Sigma ( B_1 )) \ =\  \Sigma (B_2 )$.
\par
For a canonical subalgebra $A$ of an AF $C^{\STAR} $-algebra
note that
\[
K_0 (A) \ =\lim_\to \, K_0 (A_n ) \ =
\lim_\to
\, K_0 (A_n \ \ca\  A_n^{\STAR} ) \ =\	K_0 (A \ \ca\  A^{\STAR} ).
\]
 \par
Define the
\it algebraic order
\rm $S \ =\  S(A)$ on
$ \Sigma (A)$ to be the reflexive transitive antisymmetric relation
such that $[p]S[q]$ if and only if $q\ \,\ap_{{v}}^{{}}\ p$
for some partial isometry $v$ in some algebra $A\ \,\cx\ \,M_n$ for some $n$.
For canonical subalgebras we can take $p,\ q,\ v$ in $A$, because
$K_0^+ (A)$ has cancellation.
The pair $( \Sigma (A), S(A))$ does not form
a complete invariant for such
subalgebras of AF $C^{\STAR} $-algebras, but we
shall see that it is complete for certain subclasses.
\par\noindent
\par\vspace{1.0\baselineskip}

{\noindent \bf 2. Embeddings and Normalisers.}
\par\vspace{1.0\baselineskip}
\par
Let $C \ \ib\  A \ \ib B$ be as in the
second paragraph of section 1. Every self-adjoint projection in $A$
is equivalent in $A \ \ca\  A^{\STAR} $ to a projection in
$A_n \ \ca\  A_n^{\STAR} $ for some $n$,
and so is equivalent to a projection in $C_n $ for some
$n$. Furthermore the algebraic orders can be understood in
terms of the partial isometries which normalise $C$, as
we now indicate.
\par\vspace{1.0\baselineskip}

{\noindent \bf Definition 2.1.}
The
\it normaliser
\rm of $C$ in $A$ is the semigroup
$N_C (A)$ of partial isometries $v$ in $A$ such that
$vCv^{\STAR} \ \subseteq	C$ and $v^{\STAR} Cv \ \subseteq	 C$.
The %
\it strong normaliser %
\rm of $C$ in $A$ is the
subsemigroup $N_C^s (A)$ of elements $v$
which preserve the relation $ \ap_A $ in the sense that if
$p_1 \ap_A  p_1 ,$ with $p_2 \leq v^{\STAR} v$,
and $p_2 \leq v^{\STAR} v$, then
$vp_1 v^{\STAR} \ap_A vp_2 v^{\STAR} $, and if
$p_1\ \ap_{{A}}^{{}} \ p_2 ,$ with $ \ p_1\ \le	v\ v^*\ ,
$ and $ \ p_2 \le \ v\ v^{\STAR} $,
then $v^*\ p_1\ v\ \ap _A
\ v^*\ p_2\ v$.
\par
The normaliser $N_C (A)$ has the following important property.	Each
$v$ in $N_C (A)$ has the form $cw$ with $c$ a partial isometry
in $C$ and $w$ an element of $N_{C_{k}}(A_k )$ for some $k$.
(Also, every operator of this form is in $N_C (A)$.)
We use this below without further explanation.	For details see
\cite{p-tensor} or \cite{p-book}.
\par\vspace{1.0\baselineskip}

{\noindent \bf Lemma 2.2}
Let $p$ and $q$ be projections in $A$.	Then $[p]\ S(A)\ [q]$ if
and only if there exist projections $p ^\prime$ and $q ^\prime$
in $C$ and a partial isometry $v$ in $N_C (A)$ such
that $p\  \sim \ p ^\prime \ ,\ q\  \sim \ q ^\prime$ and
$v^{\STAR} v\ =\ q\ ,\ vv^* \ =\ p$.
\par\vspace{1.0\baselineskip}

\begin{mproof}
Suppose that $[p]\ S(A)\ [q]$.	For some large $k$ there are projections
in $A_k$ which are close to $p$ and $q$ , and so it
follows that there exist projections $p ^\prime$ and $q ^\prime$
in $C_k$ with $p\  \sim \ p ^\prime \ ,\ q\  \sim \ q ^\prime$.
By the hypothesis it follows that there is a partial isometry $w$
in $A$ with $w^{\STAR} w\ =\ q ^\prime$ and $ww^* \ =\ p ^\prime$.
Increasing $k$ if necessary, choose an operator $x$ in $A_k$, close
to $w$, with $x\ =\ p ^\prime  xq ^\prime$, such that $x$ is invertible
when viewed as an operator from $q ^\prime H$ to $p ^\prime H$ where
$H$ is the finite-dimensional Hilbert space underlying $A_k$.
Let
$p ^\prime \ =\ p ^\prime_1\ \ +\ \ ...\ +\ \ p ^\prime_r\
\ ,\ \ q ^\prime\ =\ q ^\prime_1\ \ +\ \ ...\ +\ \ q ^\prime_r$ be
the decomposition into minimal
projections of $C_k$.  By the invertibility of $x$ it follows
that there is a permutation $\pi$ of $\ 1,\ ...,\ r\ $ such that
$p ^\prime_{i}\ x\ q ^\prime_{\pi (i)}$ is nonzero for each $i$.  By the
minimality of the $p ^\prime_i$ and $q ^\prime_j$ it follows that
there is a partial isometry $v_i$ with initial projection
$q_{\pi (i)}$ and final projection $p_i$.  The partial isometry
$v\ =\ v_1\ +\ ...\ +\ v_r$ satisfies
$v^{\STAR} v\ =\ q\ ,\ vv^*\ =\ p$, and since $v$ belongs to
$N_{C_{k}}(A_k )$ it follows that $v$ belongs to
$N_C (A)$.
\end{mproof}

\par
In the next section
we consider embeddings affiliated with maximal abelian
subalgebras and the following terminology
will be useful.
\par\vspace{1.0\baselineskip}

{\noindent \bf Definition 2.3.}
Let $A$, $A { \, ^\prime } $ be FDCSL algebras containing
masas $C$, $C { \, ^\prime } $ respectively and
let $ \alpha \,:\, A \rarr A { \, ^\prime } $
be an injective algebraic homomorphism with $C \rarr C { \, ^\prime } $.
The embedding is said to be
\begin{itemize}
\item[{(i)}]
\it star-extendible%
\rm , if there is an extension
$C^{\STAR} (A) \rarr C^{\STAR} (A { \, ^\prime } )$,
\item[{(ii)}]
\it regular %
\rm (with respect to $C$, $C { \, ^\prime } $) if
$N_C (A) \rarr N_{ C { \, ^\prime } } (A { \, ^\prime } )$,
\item[{(iii)}]
\it strongly regular %
\rm (with respect to
$C$, $C { \, ^\prime } $) if
$N_C^s (A) \rarr N_{{C} { \, ^\prime }}^s (A { \, ^\prime } )$,
\par\vspace{1.0\baselineskip}
\end{itemize}\par
If $A \ \subseteq   M_n $,
$A { \, ^\prime } \ \subseteq   M_n \ \cx\  M_m $,
identified with $M_n (M_m )$, and if
$A \ \cx\  \IC I \ \subseteq   A { \, ^\prime } $,
then we refer to the natural embedding
$ \rho \,:\, A \rarr A { \, ^\prime } $,
given by $ \rho (a) \ =\  a \ \cx\  1$, as a
refinement
 embedding. In particular, viewing
$T_{nm} $ as the upper triangular matrix
subalgebra of $M_n (M_m )$ we have
the {\it refinement} embedding $ \rho \,:\, T_n \rarr T_{nm} $.
In contrast, if $ \IC I \ \cx\	A \ \subseteq   A { \, ^\prime } $, then
we refer to the embedding
$ \sigma \,:\, A \rarr A { \, ^\prime } $,
given by
$ \sigma (a) \ =\  1 \ \cx\  a \ =\  a \ \oplus \ \it\hbox{...} \oplus \  a$
($n$ times) as a standard embedding.
In particular, with the same identification
$T_{nm} \ \subseteq   M_{n}( M_{m}) $, we have the
{\it standard} embeddings $ \sigma \,:\, T_m \rarr T_{nm} $.
Standard embeddings, refinement embeddings, and many
other hybrid embeddings are strongly regular.  In contrast the
embedding $T_2 \rarr T_4 $ given by

\[
\left[
\begin{array}{ll}
a & b\\
  & c
\end{array}
\right]
\quad
\to
\quad
\left[
\begin{array}{rrrr}

a & 0 & 0 & b\\
 & a & b & 0\\
 &  & c & 0\\
 &  &  & c
\end{array}
\right]
\]
is a regular star-extendible embedding which is not
strongly regular.
\par
It is straightforward to check that a strongly regular star-extendible
embedding
$T_n \rarr T_m $ is determined, up to conjugation
by a unitary in the diagonal algebra $D_m \ =\	T_m \ca (T_m )^{\STAR} $,
by its restriction to the diagonal algebra $D_n $. Recall that a
finite-dimensional nest algebra $A$ is an FDCSL algebra whose lattice of
invariant projections is totally ordered.  Similarly it can be shown
that a strongly regular embedding
$ \alpha \,:\, A \rarr A { \, ^\prime } $ between
direct sums of such algebras is
determined, up to conjugacy by a unitary in
$A { \, ^\prime } \ca (A { \, ^\prime } )^{\STAR} $, by its restriction to
$A \ca A^{\STAR} $.
\par\vspace{1.0\baselineskip}
{\noindent \bf Definition 2.4}
Let $C\ \ib\ A\ \ib\ B$ be as in the second paragraph of Section 1.  The
\it strong algebraic order
\rm $S_1 (A)$ is the subrelation of
the algebraic order $S(A)$ such that $[p] S_1 (A)\,[q]$ if and
only if there are representatives $p,\ q$ in Proj C and a
partial isometry $v$ in $N_C^s (A)$ with $q\ \ap
_{{v}}^{{}} \ p$.
\par\vspace{1.0\baselineskip}
\par
In general $S_1 (A)$ is a proper subrelation of $S(A)$.  This can be
seen for elementary FDCSL algebras. On the other hand these relations agree
in the case of triangular nest algebras. Note that $S_1 (A)$ depends
on an implicit choice of canonical masa, so it is not clear, a priori,
whether $S_1 (A)$ is even an invariant for isometric isomorphism.
However in the following triangular  context we have :
\par\vspace{1.0\baselineskip}

{\noindent \bf Lemma 2.5}
Let $A$ be the limit of the system $A_1\  \ap
\  A_2\  \ap
\  ...$
consisting of direct sums of triangular finite-dimensional nest algebras and
strongly regular embeddings.  Then $S_1 (A)\ =\ S(A)$, where $S_1 (A)$
is the strong algebraic order of $A$.
\par\vspace{1.0\baselineskip}

\begin{mproof}
If $A_k$ is a triangular finite-dimensional nest algebra then
$S_1 (A_k )\ =\ S(A_k )$.  Indeed, if $p,\ q$ are
projections in a masa $C_k$ of $A_k$, and $q\ \ap
\ p$, then we can
order the minimal subprojections (in $C_k$) of $q$ and $p$ and obtain
$q\  \ap
_{{v}}^{{}} \  p$ where  $v\ \in\ N_{C_{k}}(A_k )$
is a partial isometry which matches these subprojections in order.  Since
$v$ preserves the partial ordering on minimal projections (induced by
$A_k$) $v$ belongs to $N_{{C}_k}^s\ (A_k )$.
\par
If $[e]\ S(A)\ [f]$, and if $C$ is the limit of the subsystem
$C_1\ \ap
\ \ C_2\ \ap
\  ...$, choose $p,\ q$ in Proj$C_k$ for
some large $k$ with $[e]\ =\ [p],\  [f]\ =\ [q]$ and with
$q\ \ap
_{{A}_k}^{{}} \ p$.  Choose $v$ as above
in $N_{{C}_k}^s\ (A_k )$.  The embeddings
are strongly regular and so it follows that $v$ is in $N_C^s (A)\ $.
\end{mproof}
 \par\noindent
\par\vspace{1.0\baselineskip}


{\noindent \bf 3. Classifying limit algebras.}
\par\vspace{1.0\baselineskip}
\par
Let $A_1 \rarr A_2 \rarr \ \it\hbox{...} $ be a direct system of
FDCSL algebras with star-extendible injective embeddings.  We say
that the system has the %
\it conjugacy property %
\rm if whenever
$ \alpha \,:\, A_n \rarr A_{{n+k}} $ is a star-extendible
embedding whose restriction $ \alpha |\, C_n $ to a masa
in $A_n $ is equal to the given injection
$i \,:\, C_n \rarr A_{{n+k}} $, then there is a
unitary operator $u$ in
$A_{{n+k}} \ \ca\  A_{{n+k}}^{\STAR} $ such
that $ \alpha \ =\  (Ad\ u) \cir i $, where $i$ is the
given injection of $A_n $.  Similarly we say that the
system has the %
\it conjugacy property for strongly regular maps %
\rm if the same conclusion holds just for maps
$ \alpha \,:\, A_n \rarr A_{{n+k}} $ which are additionally
strongly regular relative to two masas.
This is the appropriate concept for direct systems with strongly
regular embeddings.
\par
We noted above that a strongly regular direct system
$A_1 \rarr A_2 \rarr \ \it\hbox{...},$ with each $A_k $ a
direct sum of triangular finite-dimensional nest
algebras, has the conjugacy property
for strongly regular maps.  On the other hand it can be shown that
the natural direct system
\[
A_1 \rarr A_1 \ \otimes\  A_2 \rarr A_1 \ \otimes\  A_2 \ \otimes\  A_3
\rarr \ \it\hbox{...},
\]
where each $A_k $ is an FDCSL algebra, has the ordinary conjugacy property.
\par
In the theorem below we obtain a generalisation of Elliott's
classification of AF
$C^{\STAR} $-algebras, and the proof is modelled on the self-adjoint case.
However, in our generality it is necessary to do extra work to lift
relations on the scale of $K_0 (A)$ to normalising partial
isometries in such a way that we obtain
star-extendible embeddings.	(We remark that even
strongly regular embeddings of FDCSL
algebras need not be star-extendible.)
\par
Let $A$ be a canonical subalgebra with canonical masa $C$ as before,
 and let
$R\ \subseteq  S(A)$ be a connected transitive reflexive finite subrelation.
We view this as a binary relation on the set $ \{\ 1 ,\ ...,\ n\ \}$.
It follows from Lemma 2.2 that we can find orthogonal projections
$ p_{1} ,\  ... ,\ p_n $ in $M_{\infty} (C)$ and partial isometries
$v_{ij}$ in $M_{\infty} (N_C (A) )$, with
$p_j\ \ap
_{{{v}_{ij}}}^{{}} \ \ p_i$, whenever
$(i,\ j)\ \in\ R$.
We say that $S(A)$ has the
\it star
realisation property
\rm if for every such subrelation there is a
choice with $ \ \{v_{ij} \ :\ (i,\ j)\	\in \ R\} \ $ a subset of a
complete matrix unit system.  By this we mean that the natural map
$A(R)\ \rightarrow \ M_{\infty} (A)$ given
by $(a_{ij} )\ \rightarrow \ (a_{ij}v_{ij} )$
is a star-extendible injection from the FDCSL algebra $A(R)$ associated
with $R$.  The star realisation property is rather restricted, as we
observe below in  Remark 3.3.  Nevertheless it
holds in several contexts of interest
and allows for the statement of a general theorem.

In an exactly
analogous way one can define when $S_1 (A)$ has the star realisation
property, and there is a corresponding variant of the following theorem
for strongly regular systems.
\par\vspace{1.0\baselineskip}

{\noindent \bf Theorem 3.1}.
Let $A$  and $A^\prime$ be the limits of the systems
$A_1 \rarr A_2 \rarr \ \it\hbox{...} $
and $A_1^\prime \rarr A_2^\prime \rarr \ \it\hbox{...} $
consisting of
FDCSL algebras and injective star-extendible regular embeddings.
Suppose further that the systems have the conjugacy property
and that  the algebraic orders $S(A)$ and $S(A ^\prime )$ have the
star realisation property. Then $A$ and
$A { \, ^\prime } $ are isometrically isomorphic if and only if there
is a scaled order group isomorphism
$ \theta \,:\, K_0 (A) \rarr K_0 (A { \, ^\prime } )$
which gives an isomorphism of the algebraic order.
\par

\begin{mproof}
Assume that $ \theta $ exists.	Thus it is assumed that $ \theta $
gives a bijection between the algebraic
orders $S(A),\ S(A ^\prime )$ associated with
canonical masas $C,\ C ^\prime$, respectively,
affiliated with the given direct systems.  It will be
enough to construct a system of embeddings
$$
A_1 \  \rightarrow
_{ \phi_1 } \  A { \, ^\prime }_{{n}_1 } \  \rightarrow
_{ \psi_1 } \  A_{{m}_1 } \  \rightarrow
_{ \phi_2 } \ \it\hbox{...}
$$
which commute with the given embeddings $A_1 \rarr A_{{m}_1 }$,
$A { \, ^\prime }_{{n}_1 } \rarr A { \, ^\prime }_{{n}_2 } ,\ \it
\hbox{...} \  $.  In particular
the constructed isomorphism
maps $\  \bigcup A_n $ onto $\  \bigcup A { \, ^\prime }_n $.
This isomorphism $ \phi \ =\  {\displaystyle \lim_\to} \phi_k $ will also
implement the given isomorphism $ \theta $.
\par
The algebra $A$ is a regular canonical subalgebra of the AF $C^{\STAR}
$-algebra
$C^{\STAR} (A) \ =\  {\displaystyle \lim_\to} C^{\STAR} (A_n )$,
and from this it follows
that we can choose systems $ \ \{e_{ij}^n \} \ $ of matrix
units for $C^{\STAR} (A_n )$, for
$n = 1,2,...$, such that each $e_{ij}^n $ is
a sum of matrix units in $ \ \{e_{ij}^{n+1} \} \ $, and such that
$ \ \{e_{ii}^n\} \ $ spans the masa $C_n $ in $A_n $.
(See \cite{p-book} for example.)
Let $ \ \{e_{ij}^n \,:\, (i,j) \memof \Omega_n\ \} $
be the set of matrix units in $A_n $.
Similarly choose the matrix unit system $ \ \{f_{ij}^n\} \ $ for
$C^{\STAR} (A { \, ^\prime }_n )$, with $ \ \{f_{ii}^n\} \ $
spanning $C { \, ^\prime }_n $, and let
$ \ \{f_{ij}^n \,:\, (i,j) \memof \Omega { \, ^\prime }_n \ \} $ be the
set $A { \, ^\prime }_n \ \ca\	\ \{f_{ij}^n\} \ $.
\par
Choose $n_1 $ large enough so that there are
orthogonal projections $g_{ii} $ in $C { \, ^\prime }_{ n_1 } $
such that $ \theta ([e_{ii}^1 ]) \ =\  [g_{ii} ]$ for all $i$.
(We need not be precise about the range of $i$.)
Let $(i,\ j)\ \in\ \Omega_1$.  Since
$ \theta $ preserves the  algebraic order,
and since $S(A ^\prime )$ has the star realisation property, there
is a choice of partial isometries $v_{ij}$ in $M_{\infty} (N_{{C} ^\prime } (A
^\prime ))$, for $(i,\ j)\ \in\ \Omega_1$, such that $g_{jj}\ \rightarrow
_{{{v}_{ij}}}^{{}} \ g_{ii}$, and such that
the induced map $\phi_1\ :\ A_1\ \rightarrow
\ M_{\infty} (A ^\prime )$ is
star-extendible.  Since the orthogonal projections $g_{ii}$ lie
in $C_{n_1}^{^\prime}$, the range of $\phi_1$ is
actually in $A ^\prime$.  In view of the remarks preceeding Lemma 2.2,
$\phi_1$ has the form
$$
\phi_1 \left( \sum_{{{\Omega}_1}}^{{}} \ a_{ij}\ e_{ij}^1 \right)\ \ =\ \
\sum_{{{\Omega}_1}}^{{}} \ a_{ij}\ c_{ij}\ w_{ij} $$
\par\noindent
where, for some large enough $k\ ,\ w_{ij}$ is a partial isometry which
is a sum of some of the matrix units in $\{f_{ij}^k\ :\ (i,\ j)\ \in\
\Omega_k^{^\prime}\} \ $, and where $ c_{ij}\ \in\  C ^\prime$ for all
$(i,\ j)\ \in\ \Omega_1$.  However, by the star-extendibility of $\phi_1$
the set $\{\ c_{ij}\ w_{ij}\ :\ (i,\ j)\ \in\ \Omega_1 \} $
is a subset of a complete matrix unit system. From this it follows that
$\{\ w_{ij}\ :\ (i,\ j)\ \in\ \Omega_1 \} $ is necessarily a subset
of a complete matrix unit system.  It can now be shown that
$c_{ij}w_{ij}\ =\ cw_{ij}c^*$ for some partial
isometry $c$ in $C ^\prime $.  So, replacing $\phi_1$ by
$(Adc^* )\cir\phi_1$, and replacing $n_1$
by $k$, we obtain the desired map $\phi_1$.  It is
 regular because the images of the matrix units of $A_1$
lie in the normaliser of a masa.
\par

We have obtained a  regular star-extendible embedding
$ \phi_1 \,:\, A_1 \rarr A { \, ^\prime }_{{n}_1 } $ such that\\
$[ \phi_1 (e_{ii}^1 )]S(A { \, ^\prime } )[ \phi_1 (e_{jj}^1 )]$ for each
matrix unit $e_{ij}^1 $ in $ A_1 $.  We now construct
the desired
map $ \psi_1 \,:\, A { \, ^\prime }_{{n}_1 } \rarr A_{{m}_1 } $.
Choose orthogonal projections $h_{ii} $ in $C_{{m}_1 } $,
for suitably large $m_1 $, so that
$[h_{ii} ] \ =\  \theta^{-1} ([g_{ii}^{{n}_1 } ]) $ for
all $i$.  We can do this in such a way so that,
for each $i$, if $g_{ii} \ =\  \sum_{{J}_i } f_{ii}^{{n}_1 } $
then $ \sum_{{J}_i } h_{ii} $ coincides with
$i(e_{ii}^1 )$.  In other words, the choice of the
projections $h_{ii} $ determine an
injection $ \omega \,:\, C { \, ^\prime }_{{n}_1 } \rarr C_{{m}_1 } $
and we can arrange this so that $ \omega  \cir \phi_1 $ agrees with the
given injection $i \,:\, C_1 \rarr C_{{m}_1 } $.
By our earlier arguments,
increasing $m_1 $ if necessary, there is a
regular star-extendible embedding
$ \hat \omega \,:\, A { \, ^\prime }_{{n}_1 } \rarr A_{{m}_1 } $
which extends $ \omega $.
Because of the hypothesised conjugacy property there is a
unitary element $u$ in $A_{{m}_1 } \ \ca \  A_{{m}_1 }^{\STAR} $
so that
$ \*q_1 \ =\  (Adu) \cir \hat \omega $ is the desired injection from
$A { \, ^\prime }_{{n}_1 } $ to $A_{ m_1 } $,
with $ \*q_1 \cir \phi_1 \ =\  i$.  Continue to
obtain the desired system.
\end{mproof}
\par\vspace{1.0\baselineskip}

It will be noticed that the star realisation property is much stronger
than  is necessary for the proof of  Theorem 3.1.
The essential point is that any finite transitive reflexive subrelation
of $S(A ^\prime )$ (respectively $S(A)$) which is isomorphic to
the relation for $A_k$ (respectively $A_k^{^\prime} )$ for
some $k$, is star realisable.

{\noindent \bf Theorem 3.2}.
Let $A$ and $A { \, ^\prime } $ be limits of direct sums of triangular
finite-dimensional nest algebras with respect to injective strongly
regular star-extendible embeddings associated with ordered Bratteli
diagrams (as in 3.8).  Then $A$ and
$A { \, ^\prime } $ are isometrically isomorphic if and only if there is a
scaled group isomorphism $ \theta \,:\, K_0 (A) \rarr K_0 (A { \, ^\prime } )$
which preserves the algebraic order.

\begin{mproof}
The proof above applies with simplifications.
Firstly, note that the maps $\phi_1, \psi_1,\dots $ are easily defined by
specifying images for the superdiagonal matrix units (the matrix
units of the first superdiagonal).
Secondly, observe that ordered
Bratteli diagram systems have the conjugacy property. Indeed, if
$\phi : E \to F$ is an ordered
Bratteli diagram embedding between triangular
elementary algebras, and if $e$ is a matrix unit in $E$, then in each summand
of $F$ the minimal subprojections of $e^*e$  interlace those of
$ee^*$.
It follows that the partial isometries $v$ in $F$ with
$e^*e \to_v ee^*$ agree modulo a multiplier of \ $C$.
\end{mproof}

{\noindent \bf Remark 3.3}
In general the algebraic order of a FDCSL algebra may not have the
star realisation property.  In fact more is true.  There are
FDCSL algebras $A_1 \ ,\ A_2$ and a scaled ordered group
injection $\theta \ :\ K_0(A_1 )\ \,\rightarrow
\ \,K_0(A_2 )$ with
$\theta^{(2)} (R(A_1 ))\ \subseteq  R(A_2 )$ which is not
induced by any regular injective embedding.  To see this consider
the subalgebra $A_2$ of matrices $(a_{ij} )$ in
$M_4 \otimes  M_2$ of the form
\[
(a_{ij}) =
\left[
\begin{array}{rrrrrrrr}
a_{11} & 0 & a_{13} & 0 & a_{15} & 0 & a_{17} & a_{18}\\
       & a_{22} & 0 & a_{24} & 0 & a_{26} & a_{27} & a_{28}\\
       && a_{33} & 0 & 0 & 0 & a_{37} & 0\\
      &&& a_{44} & 0 & 0 & 0 & a_{48}\\
     &&&& a_{55} & 0 & 0 & a_{58}\\
    &&&&& a_{66} & a_{67} & 0\\
   &&&&&& a_{77} & 0\\
  &&&&&&& a_{88}
\end{array}
\right]
\]

\par\noindent
\par\vspace{1.0\baselineskip}
\noindent Let $A_1\ =\ T_2\ \,\cx\ \,T_2$ and consider the injection
$\phi\ :\ A_1\ \ca \ A_1^* \ \,\rightarrow
\ \,A_2\  \ca \ A_2^*$ which is
given by $c\ \ \rightarrow
\ \ c\ \,\cx\ \,I_2$. This in turn induces
a map  \ $\theta \ :\ \Sigma\ (A_1 )\ \,\rightarrow
\ \,\Sigma\ (A_2 )$ with
$\theta^{(2)} \ (R(A_1 ))\ \ \subseteq  \ R(A_2 )$.
Examination shows that there is no unital regular injection
$\ A_{1}\ \rightarrow
\  A_{2}$, star-extendible or otherwise,
which induces $\theta$.  We remark that in the case of infinite
tensor products of proper finite-dimensional nest algebras the
algebraic order contains arbitrary subrelations, and in
particular subrelations isomorphic to $\theta^{(2)}\ (R(A_1 ))$.
Thus it does not seem that the methods of Theorem 3.1 are
immediately applicable in the classification of infinite tensor products.
\par\vspace{1.0\baselineskip}

{\noindent \bf Remark 3.4.}
If $A$ and $A { \, ^\prime } $ are
isomorphic as Banach algebras then
it can be shown that there is a scaled ordered group isomorphism from
$K_0 (A)$ to $K_0 (A { \, ^\prime } )$ which
preserves the algebraic order.	As a consequence
the algebraically ordered scaled ordered group $K_0 (A)$ is
a complete invariant for bicontinuous isomorphism within the classes
considered in Theorems 3.1 and 3.2.
\par
It seems plausible that any two canonical subalgebras of an AF
$C^{\STAR} $-algebra are isometrically isomorphic if they
are algebraically isomorphic.  Settling this problem will be
a good test of the effectiveness of any future methods in the
study of subalgebras of AF $C^{\STAR} $-algebras.
\par\vspace{1.0\baselineskip}

{\noindent \bf Remark 3.5.}
There exist triangular canonical subalgebras $A$ and $A ^\prime$, with
$S(A)\ =\ S_1 (A)\ ,\mbox { with }$  $\ S(A ^\prime )\ =\ S_1 (A ^\prime )$,
and with an algebraic order preserving
isomorphism $\theta\ :\ K_0 (A)\ \rightarrow
\ K_0 (A ^\prime )$,
which are nevertheless not isometrically
isomorphic.  To see this
let $B\ =\ {\displaystyle \lim_\to} (M_{{2}^k} ,\ \rho_k )$, where
$\rho_k\ :\ M_{{2}^k}\ \rightarrow
\ M_{{2}^{k+1}}\ ,\ k\ =\ 1,\ 2,\  ...$ are
refinement embeddings, and let $B ^\prime$ be the subspace
${\displaystyle\lim_\to} (M_{{2}^k}^\prime  ,\ \rho_k )$ where $M_{{2}^k}
^\prime $ is the subspace of matrices with zero diagonal.  Furthermore, let
$D\ =\ {\displaystyle \lim_\to} (D_{{2}^k} ,\ \rho_k )$
be the canonical diagonal
subalgebra of $B$.  Adopting a little notational distortion, define
$$
A\ \ =\ \  \left[
\matrix {D \cr O}\ \ \matrix {B \cr D} \right]\ \ \ ,\ \ \
A ^\prime \ \ =\ \ \left[
\matrix {D \cr O}\ \ \matrix {B ^\prime \cr D} \right]	  \ .
$$
\par\noindent
These canonical subalgebras of $M_2\ \cx\ B$ are not isometrically
isomorphic.  This can be deduced from the fact that $A$ and $A ^\prime$
do not have topologically isomorphic fundamental relations.
See \cite{p-triangular}.  We leave the
verification of the other assertions as a simple
exercise.
\par\vspace{1.0\baselineskip}

{\noindent \bf Remark 3.6.}
Theorem 3.2 is not true if the
embedding condition is relaxed.  Let
$A \ =\  {\displaystyle\lim_\to} (T_{{2}^n } , \rho )$ be the limit
of upper triangular matrix algebras with respect to
refinement embeddings, and let
$A { \, ^\prime } \ =\	{\displaystyle\lim_\to} (T_{{2}^n } , \theta_n )$ be the
limit algebra where $ \theta_n (e_{ij}^n ) \ =\	\rho (e_{ij}^n )$
if $j < 2^n $ or
$(i,j) \ =\  (2^n ,2^n )$, and
$ \theta_n (e_{{i,2}^n }^n ) \ =\   e_{{2i,}
2^{{n+1}} -1}^{{n+1}} + e_{{2i-1,2}^{{n+1}} }^{{n+1}} $,
otherwise.
These embeddings, in which the final column of
matrix units is embedded with twisted orientation, are not strongly regular.
Despite the fact that the binary relations
$( \Sigma (A),S(A))$, $( \Sigma (A { \, ^\prime } ), S(A { \, ^\prime } ))$ are
naturally isomorphic, the algebras $A$
and $A { \, ^\prime } $ are not isomorphic.
This observation is essentially due to
Peters, Poon and Wagner \cite{ppw}.  See also \cite{p-triangular}.
\par
In \cite{ppw} a partial order $<_A$ on Proj(C) is defined by
$p\ <_A\ q$ if and only if $q\ \rightarrow
_{{v}}^{{}}\ p$
for some $v$ in $N_C (A)$.  If $A$ is triangular,
so that $A\ \ca\ A^* \ =\ C$, then $<_A$ can be identified
with $S(A)$.  It is shown in \cite{ppw}, using the invariant $<_A$,
that there are uncountably many isomorphism classes of limit
algebras of the
form ${\displaystyle \lim_\to} (T{_2^k},\ \phi_k )$.  In \cite{p+w1}
related invariants are exploited in the study of nest subalgebras.
In particular it is shown that there are uncountably many nonisomorphic
triangular nest algebras $A$ in any given UHF algebra, all having
the same trace invariant $\ \{ {\rm trace} (p)\ :\ p\ \in\  {\rm Lat} \ A \} \
$. Here \ Lat$A$\  is the projection nest in $C$ determining $A$.
\par\vspace{1.0\baselineskip}

{\noindent \bf Remark 3.7.}
Baker \cite{bak} has shown that the unital limit algebras
$ {\displaystyle \lim_\to} (T_{{n}_k } , \sigma )$, associated
with standard embeddings and the
sequences $(n_k )$, with $n_k $ dividing
$n_{{k+1}} $ for all $k$, are classified by their
enveloping UHF $C^{\STAR} $-algebras.
In \cite{ppw}, \cite{p-tensor} and \cite{p-triangular} other proofs are
given and the limit algebras
${\displaystyle \lim_\to} (T_{{n}_k } , \rho )$ are similarly classified. These
standard limit algebras are special cases of the triangular
of Theorem 3.2.
\par\vspace{1.0\baselineskip}

{\noindent \bf Remark 3.8}
The limit algebras of Theorem
3.2 are determined by
\it ordered \rm Bratteli diagrams as in the following discussion.
\par
Consider, as an illustrative example, the two stationary direct systems
$A \ =\  {\displaystyle \lim_\to} (T_{{n}_k } \ \bigoplus\  T_{{n}_{{k+1}} } ,
\theta_k )$,
$A { \, ^\prime } \ =\	{\displaystyle \lim_\to}
(T_{{n}_k } \ \bigoplus\  T_{{n}_{{k+1}} }, \psi_k )$
with the strongly regular embeddings
$ \theta_k (x \ \bigoplus\  y) \ =\  y \ \bigoplus\  (x \ \bigoplus\  y)$, and
$ \psi_k (x \ \bigoplus\  y) \ =\  y \ \bigoplus\  (y \ \bigoplus\  x)$ where
$(n_k )\ =\ (1,\ 1,\ 2,\ 3,\  ... )$ is the Fibonacci sequence.
Then $C^{\STAR} (A)$ and $C^{\STAR} (A { \, ^\prime } )$ are isomorphic,
with stationary Bratteli diagram generated by

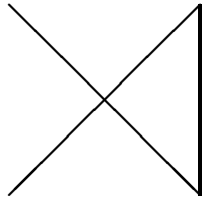
\begin{center}
\setlength{\unitlength}{0.0125in}%
\begin{picture}(80,80)(160,680)
\thicklines
\put(160,760){\line( 1,-1){ 80}}
\put(240,680){\line( 0, 1){ 80}}
\put(240,760){\line( 0, 1){  0}}
\put(240,760){\line(-1,-1){ 80}}
\end{picture}
\end{center}

However, the ordered Bratteli diagrams representing
$A$ and $A { \, ^\prime } $ are generated by the graphs

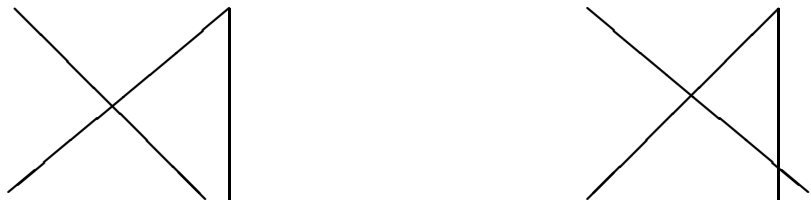
\begin{center}
\setlength{\unitlength}{0.0125in}%
\begin{picture}(330,80)(160,680)
\thicklines
\put(250,680){\line( 0, 1){ 80}}
\put(160,760){\line( 1,-1){ 80}}
\put(250,760){\line(-6,-5){ 92.459}}
\put(480,680){\line( 0, 1){ 80}}
\put(480,760){\line(-1,-1){ 80}}
\put(400,760){\line( 6,-5){ 92.459}}
\end{picture}

\end{center}

\par\noindent
and $A$ and $A { \, ^\prime } $ are not isometrically
isomorphic.  The easiest way to see this is to note that there is a
special point, $x$ say, in the Gelfand space $M(A \ \ca\  A^{\STAR} )$
with the following maximality property: if $y \ \neq \	x$
then there do not exist orthogonal projections
$ p_y $, $ p_x $ in
$A_k \ \ca\  A_k^{\STAR} $
for any $k$, such that $y( p_y ) \ =\  1$,
$x( p_x ) \ =\	1$ and
$ p_y \  \rightarrow
_{{A}_k } \  p_x $.
This point is the intersection of the supports of the
``right-most'' minimal projections in the
right summands of the $A_k $.  An isometric
isomorphism would transfer this property to
$M(A { \, ^\prime } \ \ca\  A { \, ^\prime }^{\STAR} )$, and it
is easy to check that there is no such point.
\par
The examples above fall into a class of limit algebras
associated with what might be called %
\it standard ordered Bratteli
diagrams%
\rm .  These are the Bratteli diagrams for
which at each vertex there is a
specification of the %
\it order %
\rm of the incident edges.
Such a diagram, together with a specification of the size of the
summands of $A_1$, gives rise to a unital direct system
$A_1 \rarr A_2 \rarr \ \it\hbox{...}$ in which each $A_k $
is a direct sum of upper triangular matrix
algebras. The resulting embeddings are in fact strongly regular. The examples
above illustrate the fact that these algebras are highly
dependent on the specified orderings.
\par
In a similar way one can consider ordered Bratteli diagrams
for direct sums of general finite-dimensional nest algebras.

{\noindent \bf Remark 3.9.}\ It seems quite likely that Theorem 3.2 remains
true without the assumption of triangularity. Unfortunately the proof
of this given in \cite{p-book} is incorrect because locally strongly
regular maps (ones that map matrix units into the strong normaliser)
are not necessarily strongly regular.
\par\vspace{1.0\baselineskip}

{\noindent \bf  4. Further examples}
\par\vspace{1.0\baselineskip}
\par
It is particularly interesting to calculate the algebraic orders
for nontriangular subalgebras of AF $C^{\STAR} $-algebras.
In the examples below we have a
canonical subalgebra $A$ of an AF $C^{\STAR} $-algebra $B \ =\	C^{\STAR} (A)$,
and we write $D \ =\  A \ \ca\	A^{\STAR} $ for the diagonal subalgebra. By
hypothesis, $D$ contains a canonical masa of $B$, from which it follows that
the inclusion $i \,:\, D \rarr B$ induces a surjection
$i_{\STAR} \,:\, K_0 (D) \rarr K_0 (B)$.  We shall identify this map
and the algebraic order on the scale $ \Sigma (A) \ =\  \Sigma (D)$ for
various examples.
\par\vspace{1.0\baselineskip}

{\noindent \bf  Example 4.1.}
The simplest example is finite-dimensional. Let $B \ =\  M_n $,
$A \ =\  T(n_1 ,...,n_r )$,
$D \ =\  A \ \ca\  A^{\STAR} \ =\  M_{{n}_1 }
\bigoplus \ \it\hbox{...} \bigoplus M_{{n}_r } $
where $n \ =\  n_1 + \ \it\hbox{...} + n_r $, and $A$ is the block upper
triangular subalgebra of $M_n $ associated with the ordered
$r$-tuple
$n_1 ,...,n_r $. Then   $K_0 (A) \ =\  \ZZ ^r $ with
scale $[0,n_1 ] \times \ \it\hbox{...} \times [0,n_r ]$ (with
the product order),
and
$(a_1 ,...,a_r ) S(A)(b_1 ,...,b_r )$ if and only if
$a_1 + \ \it\hbox{...} + a_r \ =\  b_1 + \ \it\hbox{...} + b_r $ and
$b_k + \ \it\hbox{...} + b_r \geq a_k + \ \it\hbox{...} + a_r $
for $1 \leq k \leq r$.	The map $i_{\STAR} \,:\, \ZZ ^r \rarr \ZZ  $ is
simply addition.  Let $ \tilde S $ be the equivalence
relation on $ \Sigma (A)$ generated by $S \ =\	S(A)$.	Then the
sets $i_{\STAR}^{-1} (x)$ for $x$ in $ \Sigma (B) \ =\	[0,n]$ are precisely the
$\tilde S $ equivalence classes.
\par\vspace{1.0\baselineskip}

{\noindent \bf  Example 4.2.}
In analogy with the last
example let $B$ be the UHF $C^{\STAR} $-algebra associated
with the generalised integer
 $2^{{\,} \infty } $, let $C$ be a canonical masa in $B$,
and consider a finite nest $0 \ < \  p_1 \ < \ \it\hbox{...} \ < p_r \ =\ $ 1
of
projections in $C$, and its associated nest subalgebra $A \ =\	\ \{b \memof B
\,:\, (1-p_j )bp_j \ =\  0, \  1 \leq j \leq r \} \ $.
Let $ \tau $ be the normalised trace on $B$ and set
$d_i \ =\  \tau (p_i - p_{{i-1}} )$,
$1 \leq i \leq r$. Then
$K_0 (B) \ =\  \IQ_d $, the binary rationals, with
the ordinary ordering, $K_0 (A) \ =\  \IQ_d^r $, and
$ \Sigma (A) \ =\  \IQ_d^r \ \ca \ ([0,d_1 ]
\times \ \it\hbox{...} \times [0,d_r ])$.
The algebraic order is exactly as in the finite-dimensional case, $i_{\STAR} $
is the addition map, and the fibres
$i_{\STAR}^{-1} (x)$ for $x$ in $ \Sigma (B)$ are
the $\tilde S $ equivalence classes.  (This latter point is a
general phenomenon.)
\par\vspace{1.0\baselineskip}

{\noindent \bf  Example 4.3.}
For a related example, let $B \ =\ {\displaystyle \lim_\to} (M_{{2}^k } , \rho
)$, let $A_k$ be the unital subalgebra
$T(n_k ,1,m_k ) \ \subseteq   M_{{2}^k } $ with
$n_k 2^{-k} \ < \  \alpha \ < \  (n_k +1) 2^{-k} $, for all $k$,
where $ \alpha $ is a fixed nondyadic point in $[0,1]$. The refinement
embeddings $ \rho $ restrict to strongly regular
embeddings $ \theta_k \ :\  A_k \rarr A_{{k+1}} $ and the limit
algebra $A \ =\  {\displaystyle \lim_\to} (A_k , \theta_k )$ can be visualised
as a subalgebra of $B$. We have
$K_0 (A) \ =\  {\displaystyle \lim_\to} K_0 (A_k ) \ =\  {\displaystyle
\lim_\to}
( \ZZ ^3 ,( \theta_k )_{\STAR} ) $ where $( \theta_k )_{\STAR} $ has the form
$$
( \theta_k )_{\STAR} \ =\  \left[ \matrix {
\matrix { 2 \cr 0 \cr 0 }
\matrix { \delta_k \cr 1 \cr 1- \delta_k }
\matrix { 0 \cr 0 \cr 2 }
}
\right]
$$
where $( \delta_k )$ is a sequence of zeros and ones.
\par
A direct argument can be given to show that $K_0 (A)$ is the subgroup of
$( \IQ_d \ +\  \alpha \ZZ  )\ \ \bigoplus\ \ ( \IQ_d \ +\  \alpha \ZZ  )$
consisting of pairs $a\ \bigoplus\ b$ with $a\ +\ b\ \,\in\ \, \IQ_d$
and that $\Sigma (A)$ is the subset with
$a\ \,\in\ [0,\,\alpha)\ ,\ b\ \,\in\ [0,\ 1\,-\,\alpha)$.
The algebraic orders agree and
$(a\  \bigoplus\  b)\ S(A)\ (c\  \bigoplus\  d)$ if and only if
$a\ +\ b\ =\ c\ +\ d$ and $b\ \ \leq  d$.
\par
One can similarly compute $( \Sigma (A), S(A))$ for
analogous algebras of the form
\[
A \ =\ {\displaystyle \lim_\to} (T(n_{{k,1}} ,...,n_{{k,r}_k } ), \rho ).
\]
\par
Notice that in all the examples above the diagonal algebra
$D$ has a certain block diagonal form. In contrast the
examples below use more interesting embeddings which
result in algebras for which $D$ is simple.  G.A. Elliott \cite{ell2}
and E.G. Effros and C.Y. Shen \cite{e+s} have analysed the dimension
groups of various stationary direct systems.  Recall that
the dimension group of a strictly positive stationary unimodular
system is  determined in terms of a distinguished
Perron-Frobenius eigenvector for the matrix determining the system.
We will not need detailed theory beyond this in the discussion below.

{\noindent \bf  Example 4.4.}
Let $A \ =\  {\displaystyle\lim_\to} (A_k , \lambda_k )$ where
$A_k \ =\  T_2 \ \cx\  M_{{4}^k } $ and
$ \lambda_k \,:\, A_k \rarr A_k \ \cx\	M_4 $
is the	regular embedding

\[
\left[
\begin{array}{ll}
x & z\\
0 & y
\end{array}
\right]
\quad
\to
\quad
\left[
\begin{array}{llll|llll}
x & & & z & 0 \\
& x & & & & z\\
&   & x & &&& z\\
&   &	& y & && & 0\\	\hline
&   &	&   & x && & z\\
&   &	&   & & y \\
&   &	&   &  & & y\\
&   &	&   &  & & & y\\
\end{array}
\right]
\]

\noindent where $x,y,z \memof M_{{4}^k } $, and where unspecified
entries are zero. Then $K_0 (A) \ =\  K_0 (D) \ =\  \IQ_d^2 $
since this is the limit of the stationary system
$$
\ZZ ^2 \rarr_{ \left[ \matrix {3 \cr 1 } \  \matrix {1 \cr 3 } \right] }
\ZZ ^2 \rarr_{ \left[ \matrix {3 \cr 1 } \
\matrix {1 \cr 3 } \right] } \ \it\hbox{...} \ \ \ .
$$
Furthermore\ \  $\IQ_d^2 $ \ has the strict ordering from the
first coordinate  ($(a,b) \leq (c,d)$ if and only if $a <  c$
or $(a,b) \ =\	(c,d)$) and the scale $ \Sigma (A)$
is the order interval $[(0,0),(1,0)]$ (see [2, page 61]).  The
map $i_{\STAR} \,:\, K_0 (A \ \ca\  A^{\STAR} )
\rarr K_0 (C^{\STAR} (A))$ is $(a,b) \rarr a+b$
($C^{\STAR} (A)$ is the
$2^{{\,} \infty }$ UHF algebra), and the algebraic
order is such that $(a,b)S(A)(c,d)$ if and only if
$a \ =\  c$ and $b \leq d$.

It is straightforward to check that the conjugacy property
holds for the system above and, more generally, for systems
over the algebras $T_2 \otimes B$ with $B$ a finite-dimensional C*-algebra.
Thus it follows from Theorem 3.1, and the remark after the proof,
that the asssociated
limit algebras are classified by the algebraically ordered scaled
$K_0$ group.

\par\vspace{1.0\baselineskip}

{\noindent \bf  Example 4.5.}
Consider the stationary system
$$
T(1,1) \ \bigoplus\  \IC \rarr_{\mu} T(1,2) \ \bigoplus\  M_2 \rarr_{\mu}
T(1,4)
\ \bigoplus\  M_3 \rarr \ \it\hbox{...} $$
where

$$
\mu \,:\, \left[
\matrix {
\matrix {x \cr 0 }
\matrix { z \cr y }
}
\right]
\ \bigoplus\  [w] \rarr
 \left[
\matrix {
\matrix { x \cr 0 \cr 0 }
\matrix { 0 \cr w \cr 0 }
\matrix { z \cr 0 \cr y }
}
\right]
 \ \bigoplus \	\left[
 \matrix {
\matrix { x \cr 0 }
\matrix { z \cr y }
}
\right]
$$
with direct limit $A \ =\   \displaystyle{\lim_\to} (A_k , \mu )$. Then
$K_0 (A) \ =\  K_0 (A \ \ca\  A^{\STAR} )$
is the limit of the stationary system

\[
{\ZZ^3} \longrightarrow _
{{\tiny
\left[
\begin{array}{ccc}
1 & 0 & 0\\
0 & 1 & 1\\
1 & 1 & 0
\end{array}
\right]}}
\quad
{\ZZ^3} \longrightarrow _
{{\tiny
\left[
\begin{array}{ccc}
1 & 0 & 0\\
0 & 1 & 1\\
1 & 1 & 0
\end{array}
\right]}}
\quad
{\ZZ^3} \longrightarrow
\dots .
\]

The embedding matrix is in $GL(3, \ZZ  )$, and so
$K_0 (A) \ =\  \ZZ ^3 $.  The enveloping
$C^{\STAR}$-algebra $B \ =\  C^{\STAR} (A)$ has
Bratteli diagram
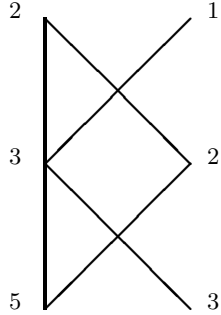
\begin{center}
\setlength{\unitlength}{0.0125in}%
\begin{picture}(89,131)(249,600)
\thicklines
\put(325,661){\line(-1,-1){ 61}}
\put(264,600){\line( 0, 1){ 61}}
\put(264,661){\line( 1,-1){ 61}}
\put(325,722){\line(-1,-1){ 61}}
\put(264,661){\line( 0, 1){ 61}}
\put(264,722){\line( 1,-1){ 61}}
\put(249,722){\makebox(0,0)[lb]{\raisebox{0pt}[0pt][0pt]{\ninrm 2}}}
\put(332,722){\makebox(0,0)[lb]{\raisebox{0pt}[0pt][0pt]{\ninrm 1}}}
\put(249,661){\makebox(0,0)[lb]{\raisebox{0pt}[0pt][0pt]{\ninrm 3}}}
\put(332,661){\makebox(0,0)[lb]{\raisebox{0pt}[0pt][0pt]{\ninrm 2}}}
\put(249,600){\makebox(0,0)[lb]{\raisebox{0pt}[0pt][0pt]{\ninrm 5}}}
\put(332,600){\makebox(0,0)[lb]{\raisebox{0pt}[0pt][0pt]{\ninrm 3}}}
\end{picture}

\end{center}

\noindent and so $K_0(B)$ is $\ZZ^2$ and
the surjection $K_0(A\cap A^\ast) \to K_0 (B)$
can be identified with the map $(\ell,m,n) \to (\ell + m,n)$. The algebraic
order is given by $(\ell,m,n) S (q,r,s)$ if and only if $n =
s, \ell + m = q + r$ and $m \leq r$.

The positive cone, $P_\alpha$ say, of $K_0(A \cap A^\ast) = \ZZ^2$,
is $\{(m,n) = m \alpha + n \geq 0\}$ where $\alpha = ( 1 + \sqrt{5})/2$,
and this, in turn,
can
be viewed geometrically as the positive
cone of the subgroup $\ZZ \alpha + \ZZ$
of $\IR$. We now indicate how to identify the positive cone ${\cal C}$ of
$\ZZ^3 =
\ZZ \oplus (\ZZ \alpha + \ZZ)$ with the set $\{(\ell, m, n) : \ell \in \ZZ_+, m
\alpha + n \in \ell (1 - \alpha) + P_\alpha\}$.

Let $p_k$ be the Fibonacci sequence so that $p_{2k+1} / p_k$
decreases to $\alpha$ and $p_k/p_{2k+1}$ increases to $\alpha$. With $Y =
\left( \begin{array}{cc} 1 & 1\\ 1 & 0 \end{array} \right)$ we have
\[
Y^k =
\left(
\begin{array}{ll}
p _{k+1} & p   _k\\
p   _k & p   _{k-1}
\end{array}
\right)
\]
\[
X^{-k} =
\left(
\begin{array}{c|cccc}
1 & & 0 & & 0\\
& \\ \hline
-p_k & & & Y^{-k}\\
p_{k+1} & & \\
\end{array}
\right)
\]

The point $(\ell,m,n)$ lies in the positive cone ${\cal C}$ if and only if for
some
$u,v,w$ in $\ZZ_+$, and for some odd integer $k$, $(\ell,m,n)^t =
X^{-k}(u,v,w)^t$. In particular, $(1,m,n)$ lies in the cone if and only if
\[
\left(
\begin{array}{c}
m \\
n
\end{array}
\right)
\quad =
\quad
\left(
\begin{array}{c}
-p_k\\
 p_{k+1}
\end{array}
\right)
+
Y^{-k}
\left(
\begin{array}{c}
v\\
w
\end{array}
\right)
\]
for some odd $k$ and $u,v$ in $\ZZ_+$. Since the smallest value of $-p_k
\alpha + p_{k+1}$ is $1 - \alpha$, and since $(0,m^\prime,n^\prime)$ lies in
$\cal C$ for all $m^\prime, n^\prime$
with $m^\prime \alpha + n^\prime \geq 0$, it
follows that $(1,m,n)$ is a point of ${\cal C}$ if and only if $m \alpha + n
\geq (1 - \alpha)$. The desired description of $\cal C$ now follows.

{\noindent \bf  Stationary Pairs of AF $C^{\STAR} $-algebras.}
\par\vspace{1.0\baselineskip}
\par
The last example is a special case of the following very general
scheme.
\par
Let $X \ =\  (a_{{ijk} \ell } )$ be an $n \times n $ matrix of
nonnegative integers, where $1 \leq k \leq k_i $,
$1 \leq \ell \leq k_j $, $1 \leq i \leq r$,
$1 \leq j \leq r$, and $k_1 + \ \it\hbox{...} +k_r \ =\  n$.
Assume that for each pair $i,j$ the partial column sum
$$
b_{ij} \ =\  a_{{ij1} \ell } + \ \it\hbox{...} + a_{{ijk}_i \ell }
$$
is independent of $ \ell $, and form the associated $ r \times r $
matrix $Y \ =\	(b_{ij} )$.  We have the commuting square of group
homomorphisms

\[
\begin{array}{ccccc}
&  & X & &\\
& \ZZ^n & \longrightarrow & &\ZZ^n\\
S & \downarrow & & S & \downarrow\\
&  & X & & \\
& \ZZ^r & \longrightarrow && \ZZ^r
\end{array}
\]

\noindent where $S$ is the homomorphism associated with the partition of $X$,
given by\\ $(S \underline x )_i \ =\  x_{ \ell_i +\,1 } + \ \it\hbox{...} + x_{
\ell_{{i+1}} } $ where
$ \ell_1 \ =\ 0$
and $ \ell_{i+1} \ =\	\ell_i\ +\ k_i$ for $i\ =\ 1, ...,\ r\,-\,1$.
For the stationary dimension
groups $G_1 \ =\  {\displaystyle\lim_\to} ( \ZZ ^n
,X)$, $G_2 \ =\  {\displaystyle\lim_\to} ( \ZZ ^r ,Y)$ we have the
induced group homomorphism $S_{\infty} \,:\, G_1 \rarr G_2 $.
Choose order units $ \underline u $ in $G_1 $ and
$\underline v \ =\  S_{\infty} \, \underline u $ in $G_2 $ and consider
AF $C^{\STAR} $-algebras $D$ and $B$ with
$K_0 (D) \ =\  G_1 $,
$K_0 (B) \ =\  G_2 $, with the chosen order units.  Furthermore
view $D$ as a unital subalgebra of $B$ so that the
inclusion map $i \,:\,D \rarr B$
induces $S_{\infty} $ ($i_{\STAR} \ =\	S_{\infty} $).
One way to visualise this inclusion is to form the stationary Bratteli
diagram for $X$ and to group together the summands associated
with the partition of the $n$ summands into $ r$ sets.
In Example 4.5 this can be indicated by the following diagram.

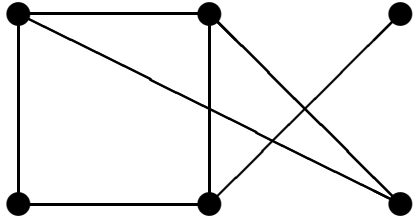
\begin{center}
\setlength{\unitlength}{0.0125in}%
\begin{picture}(170,90)(155,675)
\thicklines
\put(160,760){\circle*{10}}
\put(160,680){\circle*{10}}
\put(240,760){\circle*{10}}
\put(240,680){\circle*{10}}
\put(320,760){\circle*{10}}
\put(320,680){\circle*{10}}
\put(160,760){\line( 1, 0){ 80}}
\put(240,760){\line( 0,-1){ 80}}
\put(240,680){\line(-1, 0){ 80}}
\put(160,680){\line( 0, 1){ 80}}
\put(160,760){\line( 2,-1){160}}
\put(320,680){\line(-1, 1){ 80}}
\put(240,680){\line( 1, 1){ 80}}
\put(160,760){\line( 1, 0){ 80}}
\put(240,760){\line( 0,-1){ 80}}
\put(240,680){\line(-1, 0){ 80}}
\put(160,680){\line( 0, 1){ 80}}
\put(160,760){\line( 2,-1){160}}
\put(320,680){\line(-1, 1){ 80}}
\end{picture}
\end{center}

where the horizontal lines indicate the grouping. The partial summation
condition above is precisely the condition needed so that we can enlarge the
grouped summands to full matrix algebras, and extend the given
embeddings to these matrix algebras, and to thereby obtain a
stationary Bratteli diagram associated with $Y$.
\par
We call the resulting pair of unital AF $C^{\STAR} $-algebras
$D \ \subseteq   B$ a %
\it stationary pair%
\rm .
Actually it would be more precise to refer to the pair $D \ \subseteq   B$
as a symmetrically partitioned stationary pair, since
the same partitioning is used for rows and columns.
However we restrict attention to this symmetric case and use
the more relaxed terminology.
Clearly $D$ is a
canonical subalgebra of $B$, and so too are all the
closed intermediate algebras
$D \ \subseteq   A\subseteq   B $. It should be noted
that the pair $D \ \subseteq   B$ is determined
by the construction above, even though we make  choices
of matrix units when we form a system $B_1 \rarr B_2 \rarr \ \it\hbox{...} $
which extends the system $D_1 \rarr D_2 \rarr \ \it\hbox{...} \, $.
\par
In the example above it is easy to
see from the Bratteli diagrams that there are only two distinct proper
intermediate algebras (namely Example 4.5 and its adjoint). The
following simple pigeonhole argument shows that in
general there are only a finite number of intermediate algebras.
Let $B \ =\  {\displaystyle\lim_\to} B_n $ be the stationary unital
direct system associated with the matrix $Y$ (and a choice of order
unit) so that $D \ =\  {\displaystyle\lim_\to} D_n $ is a direct system
associated with $X$ where $D_n \ \ \subseteq   B_n $ for all $n$.
Note that for each $n$ there are exactly the same number of distinct
proper intermediate algebras, $E_n^1 ,...,E_n^s $
say, lying between $D_n $ and $B_n $. Suppose that
$D \ \subseteq   E \ \subseteq   B$. Then
$E \ =\  {\displaystyle\lim_\to} E_n $, with
$E_n \ =\  B_n \ \ca\  E$. If
$A^1 ,...,A^{{s+1}} $ are $s+1$ such algebras, then
there must exist distinct $i$ and $j$ so that
$A^i \ \ca\  B_n \ =\  A^j \ \ca\  B_n $ for
an infinity of values of $n$. Thus $A^i \ =\  A^j $.
Simple examples reveal that the number of
intermediate limit algebras
can be strictly less than $s$.
\par
We say that a stationary pair $D \ \subseteq   B$ is %
\it unimodular %
\rm if $X \memof GL(n, \ZZ  )$.  In this case $Y$
is necessarily in $GL (r, \ZZ  )$ and so both the stationary
systems for $X$ and $Y$
are unimodular in the usual sense. To see this form the matrix $X_1 $ which
is $X$ but with the rows for $i \ =\  k_1 ,k_1 + k_2 ,...,n$ replaced
by the associated partial column sums. The entries in the
new rows are the numbers $b_{ij} $ and the
determinant of $X ^\prime $ is equal to that of $X$.  Each term in the
expansion of the determinants of $X ^\prime $ along a fixed unchanged row
is divisible by $ {\rm det} \,Y$ (We can assume there is at
least one such row otherwise
$X \ =\  Y$ is either zero or has the determinant of $Y$ as a divisor.)
This can be seen more clearly if the new rows are moved by row operations to
occupy the first $r$ rows of a new matrix $ X ^{\prime\prime}  $,
still with unimodular determinant. Since $ {\rm det} \,Y$
divides $ {\rm det} \,X$ the argument is complete.
\par\vspace{1.0\baselineskip}

{\noindent \bf  Example 4.6.}
Let $X$, $Y$ be the matrices in $GL (5, \ZZ  )$ and $GL( 2, \ZZ  )$ given by

\[
X =
\left[
\begin{array}{ccccc}
1 &&&&1\\
& 1&&&1\\
&&1 && 1\\
&&&1 & 1\\
1 & 1 & 1 & 1 & 3\\
\end{array}
\right]
\quad
Y =
\left[
\begin{array}{ccc}
1 && 4\\
1 && 3
\end{array}
\right]
\]

\noindent so that $Y$ and $X$ are related as above (with $k_1 \ =\  4$,
$k_2 \ =\  1$).  Associated with the pair $X$, $Y$, and
a choice of order units, is the stationary pair
$D \ \subseteq   B$ which has joint Bratteli
diagram generated by the graph

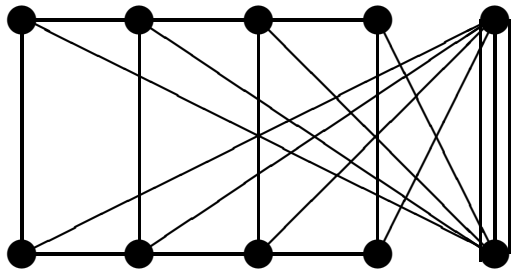
\begin{center}
\setlength{\unitlength}{0.0125in}%
\begin{picture}(210,110)(200,590)
\thicklines
\put(206,694){\circle*{12}}
\put(255,694){\circle*{12}}
\put(305,694){\circle*{12}}
\put(355,694){\circle*{12}}
\put(404,694){\circle*{12}}
\put(404,596){\circle*{12}}
\put(355,596){\circle*{12}}
\put(305,596){\circle*{12}}
\put(255,596){\circle*{12}}
\put(206,596){\circle*{12}}
\put(206,694){\line( 1, 0){149}}
\put(206,596){\line( 1, 0){149}}
\put(206,694){\line( 0,-1){ 98}}
\put(255,694){\line( 0,-1){ 98}}
\put(305,694){\line( 0,-1){ 98}}
\put(355,694){\line( 0,-1){ 98}}
\put(404,694){\line( 0,-1){ 98}}
\put(206,596){\line( 2, 1){197.600}}
\put(255,596){\line( 3, 2){148.385}}
\put(305,596){\line( 1, 1){ 98.500}}
\put(355,596){\line( 1, 2){ 47.800}}
\put(404,596){\line(-2, 1){197.600}}
\put(255,694){\line( 3,-2){148.385}}
\put(305,694){\line( 1,-1){ 98.500}}
\put(355,694){\line( 1,-2){ 49}}
\put(410,694){\line( 0,-1){ 98}}
\put(398,694){\line( 0,-1){101}}
\end{picture}

\end{center}

A choice of order units corresponds to the specification of the size
of the $5$ matrix algebra summands of $D_1 $ corresponding to
the first row of the graph. $K_0 (D) \ =\  \ZZ ^5 $ and
$K_0 (B) \ =\  \ZZ ^2 $ with positive cones $P( \underline \alpha ) $
and $P( \underline \beta ) $, respectively, determined by the eigenvectors
$ \underline \alpha \ =\  (1,1,1,1, \alpha )$,
$ \underline \beta \ =\  (1, \beta )$ for the maximal
positive eigenvalues of $X$
and $Y$. Thus
$P( \underline \alpha ) \ =\  \ \{ \underline a \memof \ZZ ^5 \,:\, (
\underline
a , \underline \alpha ) \geq  0 \ \} $,
$P( \underline \beta ) \ =\  \ \{ \underline a \memof \ZZ ^2 \,:\, \{ (
\underline
a , \underline \beta) \geq  0 \ \ \}$.
These facts follow since $Y$ is strictly positive and
$X$ has a strictly positive power. (See \cite{e+s} for more detail.)
\par
We now wish to describe all the intermediate algebras
$D \ \subseteq   A \ \subseteq   B$. The lattice
of such subalgebras is in fact a copy of the lattice of algebras
lying between $M_4 ( \IC )$ and its
diagonal subalgebra $ \IC^4 $. Indeed an algebra $E$
between $ \IC^4 $ and $M_4 ( \IC )$
is determined by a directed graph $G(E)$ on four vertices.
Fixing an assignment of these vertices to the four summands
of $D_1 $, which are grouped in $B_1 $, we can generate an
intermediate algebra $D_1 \ \subseteq   E_1 \ \subseteq   B_1 $ by
including matrix units from $B_1 $ to belong to
 $E_1  $ if there is an associated directed edge in $G(E)$. The image of $E_1$
in $B_2$ generates, with $D_2$, the analogous algebra $E_2$ and we obtain the
intermediate algebra $\tilde{E} = {\displaystyle \lim_\to} E_k$. On the other
hand if $D \subseteq A \subseteq B$ is a
closed algebra then $A \ =\  {\displaystyle\lim_\to} (A \ \ca\	B_k )$, from
which it follows that $A \ =\  \tilde E $ for some $E$. The map
$E \rarr \tilde E $ provides a bijection of intermediate algebras.
\par
Considering the special case $E \ =\  T_4 $, and some
ordering of the grouped vertices (for definiteness, take the
order corresponding to rows of $X$), we obtain an intermediate
algebra $A \ =\  \tilde E $ which is an inductive limit of
finite-dimensional nest algebras. Whilst the embeddings
$A \ \ca\  B_n \rarr A \ \ca\  B_{{n+1}} $ are {\em not} strongly
regular, they are nevertheless regular and the direct system
has the conjugacy property. Thus $A$ is an example of the algebras appearing
in Theorem 3.1. The algebraic order is given by
$ \underline a S(A) \underline b $ if and only if
$a_5 \ =\  b_5 $, $a_1 +
\ \it\hbox{...} + a_4 \ =\  b_1 + \ \it\hbox{...} + b_4 $,
and $a_i + \ \it\hbox{...} + a_4 \leq b_i + \ \it\hbox{...} + b_4 $ for
$i \ =\  1,2,3,4$.
\par
It can be shown that all of the intermediate algebras
$\tilde E \ =\ {\displaystyle \lim_\to} E_k $ are defined by systems with the
conjugacy property and have algebraic orders with the star
realisation property.  Thus they also fall within the influence of
Theorem 3.1 . As with $T_4 $, the algebraic order
of such an algebra $\tilde E $ is rather simply related to
the algebraic order of $E$.
\par
Finally it is interesting to pause to consider the fundamental
relation $R(A)$ of one of these intermediate algebras.
Recall that if $C$ is a canonical masa in $A$ then
$xR(A)y$, for $x$, $y$ in the Gelfand space $M(C)$,
if and only if there
exists $v$ in $N_C (A)$ such that $x(c) \ =\  y(v^{\STAR} cv)$
for all $c$ in $C$. It can be shown that in general $A \ =\  A^{\STAR} $
if and only if $R(A)$ is symmetric
(see \cite{th-triangular} for example). In the case of
the intermediate algebra $A \ =\  \tilde T_4 $ the equivalence
relation generated by $R(A)$ is $R(B)$.
Also the fact that there is very little room between $A$ and $B$
is reflected in the observation that $R(B)\,\setminus R(A)$
is finite. (With the exception of at most four easily identified
points $x$ in $M(C)$, the $R(A)$ orbit of $x$ and the $R(B)$
orbit of $x$ agree.)

{\noindent \bf  Final Remark.}

Recently, the interesting preprints of
Skau \cite{ska} and Herman, Putnam, and Skau \cite{hps} have appeared.
In these papers it is
shown, roughly speaking,  how ordered Bratteli diagrams provide models for
minimal homeomorphisms $\phi$ of Cantor spaces and their associated
crossed products $C(X) \times _\phi Z $. This work rests, in part,
on earlier work of Putnam \cite{put} who showed that the C*-algebra $B_x$
generated by $C(X)$ and the elements $fu$, with $u$ the canonical unitary,
and  $f$ in $C(X)$ vanishing at $x$ , is an AF C*-algebra. In
fact the tower construction in \cite{put} generates an ordered Bratteli diagram
with associated system  $ B_1 \to B_2 \to \dots $ of finite-dimensional
C*-algebras. It is not hard to check that if $A$ is the semicrossed product
\ $C(X) \times _\phi Z_+ , $\ if\ $ A_k = B_k \cap A , $ \
and if\  $ A(\phi,x) = A \cap B_x ,$
then $ A(\phi,x) = {\displaystyle\lim_\to}  A_k$
and $A(\phi,x)$ is a triangular canonical subalgebra of $B_x$ determined by a
ordered Bratteli diagram. The converse
assertion is true, and easier to establish : if $A$ is a
canonical triangular algebra determined by a  ordered Bratteli
diagram, then $A = A(\phi,x) $ for some Cantor space homeomorphism
and some point $x$ of the Cantor space.

For more on this circle of ideas see \cite{ppw2}, \cite{p-book} and the
recent preprint \cite{p+w2}.
\newpage

\vspace{1in}
e-mail \ \ maa012@central1.lancaster.ac.uk

\end{document}